\UseRawInputEncoding
\documentclass{article}
\usepackage{CJK,CJKnumb,CJKulem,times,dsfont,ifthen,mathrsfs,latexsym,amsfonts,color,indentfirst,geometry}
\usepackage{amsmath,amsthm,makeidx,fontenc,amssymb,bm,graphicx,psfrag,listings,curves,extarrows,enumitem,appendix,enumitem}
\begin{document}
%\usepackage{xeCJK}
%\setmainfont{Times New Roman}
%\setCJKmainfont{Songti SC}

%\usepackage{showkeys}
%\UseRawInputEncoding

\newtheorem{definition}{Definition}[section]
\newtheorem{theorem}{Theorem}[section]
\newtheorem*{theorem A}{Theorem A}
\newtheorem{lemma}{Lemma}[section]
\newtheorem{corollary}{Corollary}[section]
\newtheorem{proposition}{Proposition}[section]
\newtheorem{remark}{Remark}[section]

\newcommand{\al}{\alpha}
\newcommand{\ga}{\gamma}
\newcommand{\Ga}{\Gamma}
\newcommand{\dl}{\Delta}
\newcommand{\e}{\varepsilon}
\newcommand{\om}{\Omega}
\newcommand{\iy}{\infty}
\newcommand{\q}{\theta}
\newcommand{\si}{\sigma}
\newcommand{\Si}{\Sigma}
\newcommand{\la}{\lambda}
\newcommand{\vp}{\varphi}
\newcommand{\pa}{\partial}
\newcommand{\ti}{\tilde}
\newcommand{\re}[1]{(\ref{#1})}
\newcommand{\nor}[1]{\|#1\|}
\newcommand{\vs}{\vskip .25in}
\newcommand{\ov}{\overline}
\newcommand{\op}{\oplus}
\newcommand{\bs}{\backslash}
\newcommand{\zx}{\Box}
\newcommand{\ra}{\rightarrow}
\newcommand{\rh}{\rightharpoonup}
\newcommand{\hra}{\hookrightarrow}
\newcommand{\lra}{\longrightarrow}
\newcommand{\Ra}{\Rightarrow}
\newcommand{\lan}{\langle}
\newcommand{\ran}{\rangle}
\newcommand{\x}{\index}
\newcommand{\lab}{\label}
\newcommand{\f}{\frac}
\newcommand{\bt}{\begin{theorem}}
\newcommand{\et}{\end{theorem}}
\newcommand{\bl}{\begin{lemma}}
\newcommand{\el}{\end{lemma}}
\newcommand{\bd}{\begin{definition}}
\newcommand{\ed}{\end{definition}}
\newcommand{\bc}{\begin{corollary}}
\newcommand{\ec}{\end{corollary}}
\newcommand{\bprop}{\begin{proposition}}
\newcommand{\eprop}{\end{proposition}}
\newcommand{\bp}{\begin{proof}}
\newcommand{\ep}{\end{proof}}
\newcommand{\bx}{\begin{example}}
\newcommand{\ex}{\end{example}}
\newcommand{\bi}{\begin{exercise}}
\newcommand{\ei}{\end{exercise}}
\newcommand{\br}{\begin{remark}}
\newcommand{\er}{\end{remark}}
\newcommand{\be}{\begin{equation}}
\newcommand{\ee}{\end{equation}}
\newcommand{\bal}{\begin{align}}
\newcommand{\bn}{\begin{enumerate}}
\newcommand{\en}{\end{enumerate}}
\newcommand{\eal}{\end{align}}

\newcommand{\bg}{\begin{align*}}
\newcommand{\eg}{\end{align*}}
\newcommand{\bcs}{\begin{cases}}
\newcommand{\ecs}{\end{cases}}
\newcommand{\abs}[1]{\lvert#1\rvert}
\newcommand{\AR}{{\cal A}}
\newcommand{\CR}{{\cal C}}
\newcommand{\DR}{{\cal D}}
\newcommand{\FR}{{\cal F}}
\newcommand{\GR}{{\cal G}}
\newcommand{\HR}{{\cal H}}
\newcommand{\JR}{{\cal J}}
\newcommand{\LR}{{\cal L}}
\newcommand{\MR}{{\cal M}}
\newcommand{\NR}{{\cal N}}
\newcommand{\OR}{{\cal O}}
\newcommand{\PR}{{\cal P}}
\newcommand{\SR}{{\cal S}}
\newcommand{\VR}{{\cal V}}
\newcommand{\WR}{{\cal W}}
\newcommand{\A}{{\mathbb A}}
\newcommand{\F}{{\mathbb F}}
\newcommand{\N}{{\mathbb N}}
\newcommand{\X}{{\mathbb X}}
\newcommand{\Y}{{\mathbb Y}}
\newcommand{\Z}{{\mathbb Z}}
\newcommand{\C}{\mathbb C}
\newcommand{\R}{\mathbb R}
\newcommand{\RN}{\mathbb R^N}
\newcommand{\Rm}{\mathbb R^m}
\newcommand{\RNm}{\mathbb R^{N-2m}}
\newcommand{\bean}{\begin{eqnarray*}}
\newcommand{\eean}{\end{eqnarray*}}
\newcommand{\loc}{\operatorname{\rm loc}}
\newcommand{\id}{\operatorname{\rm id}}
\newcommand{\mc}{\mathcal}
\newcommand{\D}{\mathcal D}
\newcommand{\p}{\mathcal P}
\newcommand{\s}[1]{{\left( #1 \right)}}

%%%%%%%%%%%%%%%%%%%%%%%%%%%%%%%%%%%%%%%%%%%%%%%%%%%%%%%%%%%%%%%%%%%%%%%%%%%%
%%%%%%%%%%%%%%%%%%%%%%%%%%%%%%%%%%%%%%%%%%%%%%%%%%%%%%%%%%%%%%%%%%%%%%%%%%%%
%%%%%%%%%%%%%%%%%%%%%%%%%%%%%%%%%%%%%%%%%%%%%%%%%%%%%%%%%%%%%%%%%%%%%%%%%%%%

\setitemize{itemindent=38pt,leftmargin=0pt,itemsep=-0.4ex,listparindent=26pt,partopsep=0pt,parsep=0.5ex,topsep=-0.25ex}
\setlist[itemize]{leftmargin=15pt}
\numberwithin{equation}{section}

\theoremstyle{plain}

\title{\bf Normalized solutions for nonlinear Schr\"odinger systems with special mass-mixed terms: The linear couple case\thanks{This work is supported by NSFC(11801581,11025106, 11371212, 11271386),Guangdong Basic and Applied Basic Research Foundation (2021A1515010034),Guangzhou Basic and Applied Basic Research Foundation(202102020225), Province Natural Science Fund of Guangdong (2018A030310082);  E-mails: chen-z18@mails.tsinghua.edu.cn(Chen)\quad\& \quad zhongxuexiu1989@163.com(Zhong)\quad\& \quad zou-wm@mail.tsinghua.edu.cn(Zou)}}

\date{}
\author{
{\bf Zhen Chen$^1$,\;\;   Xuexiu Zhong   $^2$\;\&\;             Wenming Zou$^1$}\\
\footnotesize \it 1. Department of Mathematical Sciences, Tsinghua University, Beijing 100084, China. \\
\footnotesize \it 2. South China Research Center for Applied Mathematics and Interdisciplinary Studies, \\
\footnotesize \it  South China Normal University, Guangzhou 510631, China.}

\maketitle
\begin{center}
\begin{minipage}{120mm}
\begin{center}{\bf Abstract }\end{center}
In this paper, we prove the existence of positive solutions $(\lambda_1,\lambda_2, u,v)\in \R^2\times H^1(\R^N, \R^2)$ to the following  coupled  Schr\"odinger system
$$\begin{cases}
-\Delta u + \lambda_1 u= \mu_1|u|^{p-2}u+\beta v \quad &\hbox{in}\;\RN, \\
-\Delta v + \lambda_2 v= \mu_2|v|^{q-2}v+\beta u \quad &\hbox{in}\;\RN,
\end{cases}$$
satisfying the normalization constraints $\displaystyle\int_{\RN}u^2 =a, ~ \int_{\RN}v^2 =b$.
The parameters $\mu_1,\mu_2,\beta>0$ are prescribed  and  the masses $a,b>0$.
 Here $2+\frac{4}{N}<p,q\leq 2^*$, where $2^* = \frac{2N}{N-2} $ if $N \geq 3$ and $2^* =+ \infty $ if $N=2$. So that the terms $\mu_1|u|^{p-2}u$,$\mu_2|v|^{q-2}v$ are of the so-called mass supercritical, while the linear couple terms $\beta v, \beta u$ are of mass subcritical.  An essential novelty is that this is the first try to deal with the linear couples in the normalized solution frame with mass mixed terms, which are big nuisances due to the lack of compactness of the embedding $H^1(\R^N)\hookrightarrow L^2(\R^N)$, even working in the radial subspace. For the Sobolev subcritical case, we can obtain the existence of  positive ground state solution for any given $a,b>0$ and $\beta>0$, provided $2\leqslant N\leqslant 4$.
 For the Sobolev critical case with $N=3,4$, it can be viewed as a counterpart of the Brezis-Nirenberg critical  semilinear elliptic problem for the system case in the context of normalized solutions. Under some suitable assumptions, we obtain the existence or non-existence of  positive normalized ground state solution.
\vskip0.23in



\vskip0.23in

\end{minipage}
\end{center}

%-----------------------------------------------------------------------------------------------------------------------------------------------------------------
\vskip0.26in
\newpage
\section{Introduction}
\renewcommand{\theequation}{1.\arabic{equation}}
\renewcommand{\theremark}{1.\arabic{remark}}
\renewcommand{\thedefinition}{1.\arabic{definition}}

This paper concerns the existence of solutions $(\lambda_1, \lambda_2,u,v) \in  \mathbb R^2 \times H^1(\mathbb R^N, \mathbb R^2)$ to the following  system of elliptic equations
\be\lab{101}
\begin{cases}
-\Delta u + \lambda_1 u= \mu_1|u|^{p-2}u+\beta v \quad &\hbox{in}\;\RN, \\
-\Delta v + \lambda_2 v= \mu_2|v|^{q-2}v+\beta u \quad &\hbox{in}\;\RN, \\
\end{cases}
\ee
satisfying the normalization constraints:
\be\lab{102}
\int_{\RN}u^2 =a,   \  \  \int_{\RN}v^2 =b.
\ee
The parameters $\mu_1,\mu_2,\beta>0$ are prescribed  and  the masses $a,b>0$.
 Here $2+\frac{4}{N}<p,q\leqslant  2^*$, where $2^* = \frac{2N}{N-2} $ if $N \geqslant 3$ and $2^* =+ \infty $ if $N=2$.

\vskip0.1in

The problem under consideration is associated with the research of standing waves for the system of coupled Nonlinear Schr\"odinger Equations (NLS):
\be\lab{101.11}
\begin{cases}
-i\partial_t  \Psi_1= \Delta \Psi_1  +\mu_1|\Psi_1|^{p-2}\Psi_1 +\beta \Psi_2 \quad &\hbox{in}\; \R \times \RN, \\
-i\partial_t  \Psi_2= \Delta \Psi_2  +\mu_2|\Psi_2|^{q-2}\Psi_2 +\beta \Psi_1 \quad &\hbox{in}\;\R \times \RN. \\
\end{cases}
\ee
Systems of this type arise in nonlinear optics, for example, the propagation of optical pulses in nonlinear dual-core fiber  can be described by two linearly coupled NLS equations like
\be\lab{101.111}
\begin{cases}
i\phi_z +\phi_{xx} +|\phi|^2\phi +\kappa \psi =0,\\
i\psi_z +\psi_{xx} +|\psi|^2\psi +\kappa \phi =0,\\
\end{cases}
\ee
where $\phi(z,x)$ and $\psi(z,x)$ are the complex valued envelope functions, and $\kappa$ is the (normalized) coupling coefficient between the two cores, see \cite{ambrosetti2007multi}. Here we recall  some results    on the linearly coupled system \eqref{101} with fixed  frequency   $\lambda_i$ in the past several years. In the case of $N\leqslant 3$, $\lambda_1 =\lambda_2=\mu_1=\mu_2=1$, $p=q=4$ and $\beta >0$ small enough, Ambrosetti, Colorado and Ruiz \cite{ambrosetti2007multi} proved that \eqref{101} has multi-bump solitions. When $\mu_1|u|^{p-2}u$ and $\mu_2|v|^{q-2}v$ are replaced by $f(x,u)=(1+c(x))|u|^{p-2}u$ and $g(x,u)=(1+d(x))|v|^{q-2}v$, respsctively, system \eqref{101} has been studied by Ambrosetti \cite{ambrosetti2008remarks} with dimension $N=1$ and Ambrosetti, Colorado and Ruiz \cite{ambrosetti2008solitons} with dimension $N\geqslant2$. Brezis and Lieb \cite{brezis1984minimum} considered the more genaral case:
\be\lab{101.12}
-\Delta u_i(x)=g^i(u(x)), \ i= 1,\cdot \cdot \cdot  n,
\ee
where the $n$ functions $g^i : \RN \to  \R $ are the gradients of some function $G \in  C^1(\RN)$, namely $ g_i(u) = \partial  G(u) /  \partial u_i$. Under some conditions on $g_i$ (see (2.2), (2.3), (2.4) and (2.8) in \cite{brezis1984minimum}), they proved that \eqref{101.12} has a ground state solution which belongs to $H^1 (\RN ) \cap W^{2,q}_{loc} (\RN)$ for any $q < +\infty$ (see Theorem 2.2 and Theorem 2.3 in \cite{brezis1984minimum}). Later, Chen and Zou studied this systems with one critical exponent, see \cite{chen2012ground}.
In these works  mentioned above,    $\lambda_i$  are fixed in advance and  there is  no  constraint condition \eqref{102} imposed.

\vskip0.1in

When dealing with the Schr\"odinger equation with normalization condition in $\RN$, a new critical exponent appears, the mass-critical exponent: $\bar{p}=2+\frac{4}{N}$. Here we are going to introduce some  related works. Shibata \cite{shibata2014stable}  studied the mass-subcritical case having more general nonlinearities.
%\be \lab{101,13}
%-\Delta u + (V(x)+\lambda ) u= f(u)
%\ee
More precisely, he studied the corresponding minimizing problem with $L^2$-constraint
\be \lab{101.14}
e_a:= \inf\limits_{  u\in S_a} J(u),
\ee
where $$J(u)=\frac{1}{2}\int_{\RN}| \nabla u|^2 dx -\int_{\RN} F(|u|)dx,$$
$$S_a:= \big\{u \in H^1(\RN)  : \int_{\RN}u^2 =a\big\}$$ and $F(s)=  \int_0^s f(t)dt$ is a rather general nonlinearity. Shibata shows that there exists an $a_0 \geqslant 0$ such that $e_a$ is attained for $a > a_0$ and $e_a$ is not attained for $0 < a < a_0$.
For mass-supercritical case, the functional is unbounded from below (and above) on $S_a$. Jeanjean \cite{jeanjean1997existence} considered the problem
\be\lab{101.112}
\begin{cases}
-\Delta u +\lambda u=f(u),\\
 \int_{\RN}u^2 =a,\\
\end{cases}
\ee
with $f$ being mass supercritical and Sobolev subcritical. Under some conditions of $f$, the author obtained a radial solution at a mountain pass value. Afterward, a multiplicity result was established by Bartsch and de Valeriola in \cite{bartsch2003multi}. In the recent paper \cite{ikomatanaka2019}, Ikoma and Tanaka provided an alternative proof for this multiplicity result by exploiting an idea related to symmetric mountain pass theorems. For more general nonlinearity, see \cite{sslu2020cvpde,Mederski2020}.
Recently, the mass prescribed  problem with potential
$$\begin{cases}
-\Delta u +(V(x)+\lambda) u= f(u)\; \quad  \hbox{in}\;  \R^N,\\
 \int_{\R^N}|u|^2dx=a,
 \end{cases}$$
is also studied under some suitable assumptions on $V(x)$ by some researchers. The mass-subcritical case, we refer to Norihisa Ikoma and Yasuhito Miyamoto \cite{ikoma2020stable}, see also  Zhong and Zou \cite{zhongzou2020}. While the mass super-critical case we refer to Bartsch et al. \cite{bartsch2020}. However, they only consider the case of $f(u)=|u|^{p-2}u$.
As  for the  problem with the  combination of mass-subcritical and mass-supercritical terms, unlike the purely subcritical or supercritical situations, these mass-mixed terms make the geometric structure of the corresponding functional very complex.  Therefore  ones  need to give more detailed  arguments on the geometry of the corresponding functional constrained on suitable sub-manifold, see  \cite{soave2020,soave2020cri}.

\vskip0.12in

For the system case, Bartsch, Jeanjean, Soave \cite{bartsch2016normalized} and Bartsch, Soave \cite{bartsch2017a} considered
\be\lab{101.113}
\begin{cases}
-\Delta u +\lambda_1 u= \mu_1u^3+\beta uv^2 \quad &\hbox{in}\;\R^3, \\
-\Delta v + \lambda_2 v= \mu_2v^3+\beta u^2v \quad &\hbox{in}\;\R^3, \\
\int_{\RN}u^2 =a, \int_{\RN}v^2 =b
\end{cases}
\ee
with $\beta >0$ and $\beta <0$ respectively, they found solutions to \eqref{101.113} for specific range of $\beta$ depending on $a, b$. Recently, Bartsch, Zhong and Zou \cite{bartschzhongzou2020} studied  \eqref{101.113}  for  $\beta >0$ belongs to much  better ranges independent of the masses  $a$ and $b$, they  adopted  a new approach based on the fixed point index in cones and the bifurcation theory. For more results about \eqref{101.113}, see also \cite{bartsch2019multiple} and the references therein. When considering \eqref{101.113} with more general exponents of the following type:
\be\lab{101.114}
\begin{cases}
-\Delta u +\lambda_1 u= \mu_1|u|^{p-2}u+ \beta |u|^{r_1-2}|v|^{r_2} u \quad &\hbox{in}\;\RN, \\
-\Delta v + \lambda_2 v= \mu_2|v|^{q-2}v+\beta |u|^{r_1}|v|^{r_2-2} v \quad &\hbox{in}\;\RN, \\
\int_{\RN}u^2 =a, \int_{\RN}v^2 =b,
\end{cases}
\ee
 Gou and Jeanjean \cite{gou2016existence} proved the pre-compactness of the minimizing sequences up to translation for mass-subcritical problems. Ikoma \cite{ikomasys2020} also considered this system with potentials. As for the mass-supercritical case, we refer to \cite{bjl2018systems,lihwz2020,goujean2018}.

\vskip0.13in
\br
We emphasize that an essential novelty of this paper is that, this is the first try to dealing with the linear coupled terms in the mass-mixed situation.
Compared with the case of a single equation, there is a very complex situation in the system: when linear coupling is involved.  As far as we know, there is only one work (by Chen and Zou\cite{chen-zou2021jmaa}) concerning the linear couples in the normalized solution frame. Due to lack of compactness of the embedding $H^1(\R^N)\hookrightarrow L^2(\R^N)$, even working in the radial subspace. However, they are only focused on the mass sub-critical case, and thus the geometric structure of the functional is simple. In the present paper, we concern about \eqref{101}-\eqref{102} with $2+\frac{4}{N}<p,q\leqslant  2^*$, which means that $\mu_1|u|^{p-2}u$ and $\mu_2|v|^{q-2}v$ are of mass-supercritical terms. On the other hand, linear couple terms are of mass-subcritical. So that we are dealing with a special mass-mixed case, and the geometric structure of the functional will be very complicated.
\er

Our result is also obtained by looking for a critical point of $J$ constrained on $S(a,b)$. However, the procedure is now more complex. Passing from mass-subcritical case to supercritical case implies a striking modification in the geometry of the variational problem. Indeed, we now have
\be
\inf\limits_{ (u,v) \in S(a,b)} J(u,v) = - \infty,
\ee
thus, it is not expected to search for a minimum of $J$ on $S(a,b)$. We need to seek for a critical point having a mini-max characterization.
We introduce a Pohozaev-type constraint as follows
\be \lab{105}
{\cal P} := \{ (u,v) \in H^1(\RN) \times H^1(\RN)\setminus{\{0\}}, P(u,v)=0   \},
\ee
where
\be \lab{106}
P(u,v)=\int_{\RN}(|\nabla u|^2+ |\nabla v|^2) dx -\mu_1 \gamma_p\int_{\RN}|u|^pdx -\mu_2\gamma_q \int_{\RN} |v|^q dx.
\ee
We shall see that ${\cal P} \cap S(a,b)$ contains all solutions of \eqref{101}-\eqref{102}.
Hence we could define
\be\lab{eq:m-ab}
m(a,b) := \inf\limits_{ (u,v) \in {\cal P} \cap  S(a,b)} J(u,v) > - \infty .
\ee
We note that $m(a,b)$  is actually a mini-max value of $J(u,v)$ on $S(a,b)$. However, it is very difficult to check the complete mini-max structure for  arbitrary $\beta >0$, since the presence of the linear couple terms (It seems that one can check the mountain pass geometry structure provided $\beta>0$ small enough).

\vskip0.12in

On the other hand, as mentioned above, it is difficult to verify the compactness of  $\{(u_n, v_n)\}$. We need to do much more complex compactness analysis and energy estimation, due to the linear coupled terms. However, if we define
\be\lab{eq:m-radial}
m_r(a,b):=\inf\limits_{ (u,v) \in {\cal P} \cap  S_r(a,b)} J(u,v),
\ee
we can prove that $m_r(a,b)=m(a,b)$ (see Lemma \ref{lemma:radial} below).
In  the  presence paper, we are aim to prove the achievement of $m(a,b)$.  So we can work in the framework of radial subspace.

\vskip0.12in

 In order to gain a P-S sequence of $J$ at the level $m(a,b)$, we apply the well known Caristi fixed point theorem to construct one. Precisely,  by the fiber map $u \to  t\star u=  t^{\frac{N}{2}}u(tx)$ for $(t,u) \in \R^+ \times H^1(\RN)$, we could construct a $C^1$ map from $S_r(a,b)$ to  $ {\cal P} \cap  S_r(a,b)$:
\be
(u,v) \to t_{(u,v)} \star (u,v).
\ee
Define $\varphi(u,v) :=J(t_{(u,v)} \star (u,v))$, by Caristi fixed point theorem  and Ekland variational principle,   we could find  a P-S sequence ${(u_n, v_n)}$ of $J(u, v)$ at the level $m(a, b)$, with the additional conditions $(u_n, v_n) \in {\cal P}$, see Lemma\ref{l301} and Corollary \ref{c301} in Section \ref{sec:Existence-PS-sequence}.

\vskip0.12in

Here comes our first main theorem, focused on the case $2+\frac{4}{N}<p,q<2^*$ and $2\leqslant N\leqslant 4$.

\bt \lab{thm1}
Suppose  $2 \leqslant N \leqslant 4, \mu_1,\mu_2>0$ and $2+\frac{4}{N}<p,q< 2^*$. Then for any $a>0,b>0,\beta>0$,  problem \eqref{101} has a  ground state  solution $(\bar{\lambda}_1,\bar{\lambda}_2, \bar{u},\bar{v})$ with  $\bar{\lambda}_1, \bar{\lambda}_2 >0$ and  $\bar{u},\bar{v}$ are both positive and radial, satisfying the constraint \eqref{102}. That is,
\be
\begin{split}
J(\bar{u},\bar{v}) &=\inf \{ (u,v): (u,v) \in {\cal P} \cap S(a,b)    \} \\
& =\inf \{ (u,v): (u,v) \ \  \text{is a solution of \eqref{101}-\eqref{102} for some }\lambda_1,\lambda_2   \}
\end{split}
\ee
holds.
\et

\br
If $\beta <0$, then \eqref{101}-\eqref{102} doesn't have positive  ground state  solutions. Indeed, if $\bar{u},\bar{v} >0$ solves \eqref{101}-\eqref{102} with $J(\bar{u},\bar{v})=m(a,b)$, then $J(-\bar{u},\bar{v})=J(\bar{u},-\bar{v})< J(\bar{u},\bar{v})=m(a,b)$,  which is a contradiction with $m(a,b) \leqslant J(-\bar{u},\bar{v})=J(\bar{u},-\bar{v}) $, due to the observation $P(-\bar{u},\bar{v})=P(\bar{u},-\bar{v})= P(\bar{u},\bar{v})=0$.
\er

If one of the exponent $p,q$ is Sobolev critical, things become tougher.  Let  $\gamma_p=\frac{N(p-2)}{2p}$ and denote ${\cal C}_{N,p}, { \cal S}$ be the best constants in
Gagliardo-Nirenberg inequality and Sobolev inequality respectively,   we have the following results.

\bt \lab{thm2}
Suppose  $ N=3,4$  and assume $  2+\frac{4}{N}<p< 2^*, q =2^* $. Let $\beta, \mu_1, \mu_2>0$  be fixed. If $a>0,b>0$ satisfies

\be\lab{zwm=1}  m(a,b)+\beta \sqrt{ab}<\frac{1}{N}\frac{{ \cal S}^{\frac{N}{2}}}{\mu_2^{\frac{N}{2}-1}},\ee
 where $m(a,b)$ is defined by \eqref{eq:m-ab}, then
 problem \eqref{101} has a  ground state  solution $(\bar{\lambda}_1,\bar{\lambda}_2, \bar{u},\bar{v})$ with  $\bar{\lambda}_1, \bar{\lambda}_2 >0$ and  $\bar{u},\bar{v}$ are both positive and radial, satisfying the constraint \eqref{102}. That is,
\be
\begin{split}
J(\bar{u},\bar{v}) &=\inf \{ (u,v): (u,v) \in {\cal P} \cap S(a,b)    \} \\
& =\inf \{ (u,v): (u,v) \ \  \text{is a solution of \eqref{101}-\eqref{102} for some }\lambda_1,\lambda_2   \}
\end{split}
\ee
holds. Moreover, $\bar{\lambda}_1, \bar{\lambda}_2 >0$, and  $\bar{u},\bar{v}$ are both positive and radial.
\et

The condition (\ref{zwm=1}) looks not so  natural,  however we can prove the monotonicity of $m(a,b)$ with respect to $b$, consequently, we have the following result.
\begin{corollary}\lab{cro-th2}
Suppose  $ N=3,4$  and assume $  2+\frac{4}{N}<p< 2^*, q =2^* $. Let $\beta, \mu_1, \mu_2$  be fixed. If $a>0,b>0$ satisfies
\be \lab{106.2}
\Big( \frac{1}{2}-\frac{1}{p\gamma_p} \Big) \Big( \gamma_p  {\cal C}_{N,p}\mu_1 a^\frac{p-p\gamma_p}{2} \Big)^\frac{2}{2-p\gamma_p} +\beta \sqrt{ab}<  \frac{1}{N}\frac{{ \cal S}^{\frac{N}{2}}}{\mu_2^{\frac{N}{2}-1}},
\ee
then the results of  Theorem \ref{thm2} hold.
\end{corollary}

\br Corollary \ref{cro-th2} means  that  the results of  Theorem \ref{thm2} hold when $a$ is large and $b$ is small enough.\er

Finally, for the case of $p=q =2^*$, we have the following nonexistence result.

%\br (Non existence).%
\bt \lab{thm3}(Non-existence)
Suppose  $ N =3, 4$ and assume $ \beta >0, p=q= 2^*$. Let $a,b,\mu_1, \mu_2$  be fixed, then  \eqref{101}-\eqref{102}  has no nontrivial nonnegative solutions.
\et

The rest of this paper is organized as follows. In Section \ref{sec:preliminaries}, we display some preliminary results. In Section \ref{sec:Existence-PS-sequence} and Section \ref{sec:boundedness}, we construct a bounded nonnegative radial P-S sequence $\{(u_n,v_n)\}$ for $J$ at the level $m(a,b)$ .  In Section \ref{sec:compactness}, we study the compactness of the P-S sequence $\{(u_n,v_n)\}$ given by Section \ref{sec:Existence-PS-sequence}. In Section \ref{sec:proofs}, we complete the proofs of our main theorems.

\br\lab{remark:N-2-4}
Indeed, the results established in Sections \ref{sec:preliminaries}-\ref{sec:boundedness} are valid for general dimension $N\geqslant 1$.
 However, we require $2\leqslant N\leqslant 4$ in Section \ref{sec:compactness} to do the compactness analysis, since the compact embedding for the radial space is not valid for $N=1$. On the other hand, for the technique reason, we need a Liouvill type result due to Ikoma (see \cite[Lemma A.2]{ikoma2014compactness}):
 $$-\Delta u\geqslant 0\;\hbox{in}\;\R^N, u\in L^{s}(\R^N)\;\hbox{with}\;s\begin{cases}\in (0, \frac{N}{N-2})\quad &N\geqslant3\\
 \in (0,\infty)& N=1,2 \end{cases}
 \Rightarrow u\equiv 0.$$
 By $2\leqslant \frac{N}{N-2}$, we require that $N\leqslant 4$.
\er

Throughout the paper we use the notation $|\cdot|_p$ to denote the $L^p(\RN)$ norm, and we simply write $H^1 =H^1(\RN), H = H^1(\RN) \times H^1(\RN)$. Similarly, $H_r^1$ denotes the radial subspace, and $H_r = H_r^1  \times H_r^1$,
$S_r(a,b) =S(a,b) \cap H_r$. The symbol $|| \cdot ||$ denotes the norm in $H^1$ or $H$. Let $ u^* $ be the symmetric decreasing rearrangement of $u \in H^1$, we recall that (see \cite{lieb2001analysis})
for $1\leqslant p<\infty$,
\be \lab{113}
|\nabla u^*|_2 \leqslant |\nabla u|_2, |u^*|_p =|u|_p,\  \text{and} \ \int_{\RN}uvdx \leqslant \int_{\RN}u^*v^*dx.
\ee
The notation $\rh$ denotes weak convergence in $H^1$ or $H$. Capital latter $C$ stands for positive constant which may depend on $N,p,q\ldots$, whose precise value can change from line to line.

%-----------------------------------------------------------------------------------------------------------------------------------------------------------------
\vskip0.26in
\section{Preliminaries}\lab{sec:preliminaries}
\renewcommand{\theequation}{2.\arabic{equation}}
\renewcommand{\theremark}{2.\arabic{remark}}
\renewcommand{\thedefinition}{2.\arabic{definition}}

In this section, we summarize several results that will be used in the rest discussion. \\

Let $p \in [2,\frac{2N}{N-2})$ if $N\geqslant 3$ and $p\geqslant 2$ if $N=1,2$, let us recall the well known Gagliardo-Nirenberg inequality
\be \lab{201}
|u|_p \leqslant {\cal C}_{N,p}   |\nabla u|_2^{\gamma_p}|u|_2^{1- \gamma_p}, \ \ \forall u \in H^1,
\ee
where $\gamma_p =\frac{N(p-2)}{2p}$.

For every $N \geqslant 3$, there exists an optimal constant ${\cal S}>0$ depending only on $N$ such that
\be \lab{201.1}
{\cal S}|u|_{2^*}^2 \leqslant |\nabla u|_2^2, \ \ \forall u \in {\cal D}^{1,2}(\RN),
\ee
where ${\cal D}^{1,2}(\RN)$ denotes the completion of $C_c^{\infty}(\RN)$  to the norm $||u||_{{\cal D}^{1,2}} :=  |\nabla u|_2$. \eqref{201.1} is also named Sobolev inequality.

 Let  $a > 0, \mu > 0, 2+\frac{4}{N}<p<2^*$ be fixed,  $ (\lambda, u)\in \R\times H^1$  solves the following equation.
\be\lab{202}
\begin{cases}
-\Delta u + \lambda u= \mu |u|^{p-2}u \quad &\hbox{in}\;\RN, \\
\int_{\RN}u^2 =a.
\end{cases}
\ee
Solutions of \eqref{202} can be found as critical points of $J_{\mu,p}: H^1 \to \R$
\be \lab{203}
J_{\mu,p}(u) := \frac{1}{2}\int_{\RN} |\nabla u|^2 dx - \frac{\mu}{p}\int_{\RN} |u|^pdx
\ee
constrained on $S_a$, and the parameter $\lambda$ appears as Lagrange multiplier.  Using scaling, \eqref{202} is equivalent to
\be \lab{204}
-\Delta w + w= |w|^{p-2}w \quad \text{in}\;\RN, \quad w \in H^1, \\
\ee
whose positive solutions had been studied clearly.
So the existence of normalized solutions of \eqref{202} can be obtained easily, but there are special properties of these solutions as a normalized solution. To be precise, we introduce the Pohozaev-type constraint for single equations
\be \lab{205}
{\cal P}_{\mu,p} := \{ u \in H^1\setminus{\{0\}}: P_{\mu, p}(u) =0   \},
\ee
where
\be \lab{206}
P_{\mu, p}(u) := \int_{\RN} |\nabla u|^2 dx- \mu \gamma_p \int_{\RN}|u|^p dx.
\ee
We consider the following constraint minimization problem
\be \lab{207}
m_p^{\mu} (a) = \inf\limits_{u \in S_a \cap {\cal P}} J_{\mu,p}(u).
\ee
Then we have the following lemma.
\bl \lab{l201}
Suppose $N \geqslant 1, \mu >0,a>0$, and $ 2+\frac{4}{N}<p<2^*$, then up to a translation, \eqref{202} has an unique positive solution $u_{\mu,p,a} $ with $\lambda >0$. Moreover,
$$m_p^{\mu} (a)= J_{\mu,p}(u_{\mu,p,a} ) >0,$$
and $m_p^{\mu}(a)$ is strictly decreasing with respect to $a$.
\el
\bp
By \cite{gidas1979symmetry,kwong1989uniqueness}, up to a translation, $w$ is the unique positive solution of (2.7), which is radially symmetric and decreasing with respect to 0. Since $2+ \frac{4}{N}<p<2^*$, by scaling, we obtain the unique solution of \eqref{202}
$$u_{\mu,p,a}=\Big(\frac{\lambda}{\mu}\Big)^{\frac{1}{p-2}}w (\lambda^{\frac{1}{2}}x)  \ \ \text{with} \ \  \lambda= \Big(\frac{a}{|w|_2^2}\mu^{\frac{2}{p-2}}\Big)^\frac{p-2}{2-p\gamma_p}.$$     \\
Then
\be \lab{206.1}
\begin{aligned}
m_p^{\mu} (a) &= J_{\mu,p}(u_{\mu,p,a} )\\
&=\Big( \frac{1}{2}-\frac{1}{p\gamma_p} \Big) \int_{\RN} |\nabla u_{\mu,p,a}|^2 dx  \\
&=\Big( \frac{1}{2}-\frac{1}{p\gamma_p} \Big) \Big( \gamma_p  {\cal C}_{N,p}\mu a^\frac{p-p\gamma_p}{2} \Big)^\frac{2}{2-p\gamma_p} \\
&>0,
\end{aligned}
\ee
and $m_p^{\mu}(a)$ is strictly decreasing about $a$ since $2+ \frac{4}{N}<p<2^*$.
\ep

 Here we recall the following lemmas without proof, due to Brezis and Lieb \cite{brezis1983a}.

\bl \lab{l202}
Suppose $(u_n, v_n) \subset H$ is a bounded sequence, $ (u_n, v_n)\rh
(u,v)$ in $ H$, then for $0<p<\infty$, we have
\be \lab{231}
\begin{aligned}
\lim\limits_{n\to \infty} \int_{\RN} |\nabla u_n|^2 -|\nabla u|^2- |\nabla (u_n-u)|^2 dx=0, \\
\lim\limits_{n\to \infty} \int_{\RN} |\nabla v_n|^2 -|\nabla v|^2- |\nabla (v_n-v)|^2 dx=0,
\end{aligned}
\ee
\be \lab{232}
\begin{aligned}
\lim\limits_{n\to \infty} \int_{\RN} |u_n|^p -|u|^p- |u_n-u|^p dx=0, \\
\lim\limits_{n\to \infty} \int_{\RN} |v_n|^p -|v|^p- |v_n-v|^p dx=0,
\end{aligned}
\ee
\be \lab{233}
\lim\limits_{n\to \infty} \int_{\RN} u_nv_n -uv- (u_n-u)(v_n-v) dx=0.
\ee
\el

\bl\lab{lemma:radial}
$\displaystyle m_r(a,b)=m(a,b)$, where $m_r(a,b)$ is given by \eqref{eq:m-radial}.
\el
\bp
It is trivial that $m_r(a,b)\geqslant m(a,b)$. For any $(u,v)\in {\cal P} \cap  S(a,b)$, we also have $(u^*, v^*)\in S_r(a,b)$.  Furthermore, $|\nabla u^*|_2^2\leqslant |\nabla u|_2^2, |\nabla v^*|_2^2\leqslant |\nabla v|_2^2$ by the P\'olya-Szeg\"o inequality,  and $\displaystyle\int_{\R^N}u^*v^*dx\geqslant \int_{\R^N}uvdx$ due to the Hardy-Littlewood inequality.
So we have
\begin{align*}
J(t\star(u,v))=&\frac{t^2}{2}\left(|\nabla u|_2^2+ |\nabla v|_2^2\right)-\frac{\mu_1}{p}t^{p\gamma_p} |u|_p^p -\frac{\mu_2}{q}t^{q\gamma_q}|v|_q^q -\int_{\R^N}uvdx\\
\geqslant&\frac{t^2}{2}\left(|\nabla u^*|_2^2+ |\nabla v^*|_2^2\right)-\frac{\mu_1}{p}t^{p\gamma_p} |u^*|_p^p -\frac{\mu_2}{q}t^{q\gamma_q}|v^*|_q^q -\int_{\R^N}u^*v^*dx\\
=&J(t\star(u^*,v^*)), \forall t>0
\end{align*}
and thus,
$$J(u,v)=\sup_{t>0}J(t\star(u,v))\geqslant \sup_{t>0}J(t\star(u^*,v^*))\geqslant m_r(a,b).$$
By the arbitrary of $(u,v)\in {\cal P} \cap  S(a,b)$, we deduce that $m(a,b)\geqslant m_r(a,b)$.
Hence, $m(a,b)=m_r(a,b)$.
\ep

\bl\lab{lemma:bounded-away}
Suppose $2+\frac{4}{N}<p,q\leq 2^*$.
$S(a,b)\cap \mathcal{P}$ is bounded away from $(0,0)$ in $D_{0}^{1,2}(\R^N)\times D_{0}^{1,2}(\R^N)$. In another word, there exists some $C>0$ such that
$$|\nabla u|_2^2+|\nabla v|_2^2\geqslant C,\;\forall (u,v)\in S(a,b)\cap {\cal P}.$$
\el
\bp
By the G-N inequality,  for $\forall (u,v)\in \mathcal{P}\cap S(a,b)$, we have
\be \lab{410}
\begin{split}
| \nabla u|_2^2 +|\nabla v|_2^2 =&\mu_1 \gamma_p |u|_p^p + \mu_2 \gamma_q |v|_q^q \\
\leqslant & C(p,q,a,b) \left(|\nabla u|_{2}^{p\gamma_p }+|\nabla v|_{2}^{q\gamma_q }\right).
\end{split}
\ee
Observing that $p\gamma_p, q\gamma_q>2$, we see that
\be \lab{411}
|\nabla u|_2^2+|\nabla v|_2^2\geqslant C >0,
\ee
\ep

 We define, for $t \in \R^+ $ and $u \in H^1 $, the radial dilation
\be \lab{301}
(t \star u)(x) := t^{\frac{N}{2}}u(tx).
\ee
It is straightforward to check that if $ u\in S(a)$, then $ (t \star u)(x) \in S(a)$ for every $t \in \R^+ $. Furthermore, the map $(t,u) \in \R \times H^1 \to t \star u \in H^1$ is continuous.

For $(u,v) \in S(a,b)$, consider the $C^1$ functional
\be \lab{302}
\begin{split}
\tilde{J}(t,u,v) &=J\s{t \star \s{u,v}} \\
&= \frac{t^2}{2}\int_{\RN}| \nabla u|^2 +|\nabla v|^2 dx  -\frac{\mu_1}{p}t^{p\gamma_p} \int_{\RN}|u|^p dx -\frac{\mu_2}{q} t^{q\gamma_q}\int_{\RN}|v|^q dx -\beta \int_{\RN}uv dx,
\end{split}
\ee
where $t \star (u,v) :=(t \star u,t \star v)$ for short. A direct computation shows that,
\be \lab{303}
\frac{d}{dt}J\big(t \star (u,v)\big)= \frac{1}{t}P\big(t \star (u,v)\big).
\ee
 Hence, for any $(u,v) \in S(a,b)$, there exists an unique $t_{(u,v)} >0$, such that $P\big(t_{(u,v)} \star (u,v)\big)= 0$ and $P\big(t \star (u,v)\big)> 0$ for $t< t_{(u,v)}$, $P\big(t \star (u,v)\big)< 0$ for $t> t_{(u,v)}$. On the other hand, applying a similar argument as \cite{soave2020},  it is not difficult to check that the set
\be
\begin{split}
{\cal P}^0:&= \{ (u,v) \in {\cal P}:  2\int_{\RN}(|\nabla u|^2+ |\nabla v|^2) dx = \mu_1 p\gamma_p^2\int_{\RN}|u|^pdx + q\mu_2\gamma_q^q \int_{\RN} |v|^q dx   \}\\
&=\{ (u,v) \in {\cal P}: \frac{d^2}{dt^2}J\big(t \star (u,v)\big)|_{t=1} =0      \}
\end{split}
\ee
is empty, then by the  implicit function theorem, the map $(u,v) \to t_{(u,v)}$ is of class ${\cal C}^1$.

%\bl\lab{lemma:bounded-away-2}
%There exists a positive number $\varepsilon_0>0$ such that
%$$0<\varepsilon_0\leqslant \inf_{(u,v)\in S(a,b)\cap {\cal P}}t_{(u^*, v^*)}\leqslant \sup_{(u,v)\in S(a,b)\cap {\cal P}}t_{(u^*, v^*)}=1.$$
%\el
%\bp
%For any $(u,v)\in S(a,b)\cap {\cal P}$, we also have that $(u^*,v^*)\in S(a,b)$. By
%\begin{align*}
%0=&P(t_{(u^*, v^*)}\star(u^*,v^*))\\
%=&t_{(u^*, v^*)}^2 \left(|\nabla u^*|_2^2+|\nabla v^*|_2^2\right)-\mu_1\gamma_pt_{(u^*, v^*)}^{p\gamma_p}|u^*|_p^p -\mu_2\gamma_q t_{(u^*, v^*)}^{q\gamma_q}|v^*|_q^q\\
%=&t_{(u^*, v^*)}^{2}\left[|\nabla u^*|_2^2+|\nabla v^*|_2^2-\mu_1\gamma_pt_{(u^*, v^*)}^{p\gamma_p-2}|u|_p^p -\mu_2\gamma_q t_{(u^*, v^*)}^{q\gamma_q-2}|v|_q^q\right]\\
%\leqslant&t_{(u^*, v^*)}^{2}\left[|\nabla u|_2^2+|\nabla v|_2^2-\mu_1\gamma_pt_{(u^*, v^*)}^{p\gamma_p-2}|u|_p^p -\mu_2\gamma_q t_{(u^*, v^*)}^{q\gamma_q-2}|v|_q^q\right],
%\end{align*}
% and
% $$|\nabla u|_2^2+|\nabla v|_2^2=\mu_1\gamma_p|u|_p^p +\mu_2\gamma_q |v|_q^q,$$
% we have that $t_{(u^*, v^*)}\leqslant 1$.
%
% then there exists an unique $t_{u^*,v^*}>0$ such that
%\ep

%-----------------------------------------------------------------------------------------------------------------------------------------------------------------
\vskip0.26in
\section{Existence of radial Palais-Smale sequence}\lab{sec:Existence-PS-sequence}
\renewcommand{\theequation}{3.\arabic{equation}}
\renewcommand{\theremark}{3.\arabic{remark}}
\renewcommand{\thedefinition}{3.\arabic{definition}}

Define a functional $ \varphi : S_r(a,b) \to \R \cup \{  + \infty\}$,
\be
 \varphi (u,v) = \max\limits_{t>0}J\s{t \star \s{u,v}}=\max\limits_{t>0}\tilde{J}(t,u,v),
\ee
one can see that $\varphi$ is continuous and bounded below.  By the Ekeland's variational principle, we have the following statement.

\bl  \lab{l301}
For $\forall \varepsilon >0$, and $(u_0,v_0) \in S_r(a,b)$ satisfies $\varphi (u_0,v_0) \leqslant \inf\limits_{(f,g) \in S_r(a,b)} \varphi(f,g) + \varepsilon$, then for any $ \delta >0$, there exists some $(u_{\delta},v_{\delta}) \in S_r(a,b),  u_{\delta} \geqslant  0, v_{\delta} \geqslant  0$, which satisfies
\be \lab{304}
\varphi(u_{\delta},v_{\delta}) \leqslant \varphi (u_0,v_0),
\ee
\be \lab{305}
|| (u_{\delta},v_{\delta}) -(u_0,v_0)||_{H} \leqslant \delta,
\ee
\be  \lab{306}
\varphi(u,v) > \varphi(u_{\delta},v_{\delta}) - \frac{\varepsilon}{\delta}|| (u_{\delta},v_{\delta}) -(u,v)||_{H}, \quad \forall (u,v)\in S_r(a,b)\backslash (u_\delta, v_\delta).
\ee

\el

\bp
Firstly, define
\be  \lab{307}
X_1=\left\{(f,g)\in S_r(a,b): f\geqslant  0, g\geqslant 0, \varphi(f,g)+\frac{\varepsilon}{\delta}\|(u_0,v_0)-(f,g)\|_H\leqslant \varphi(u_0,v_0)\right\},
\ee
then $X_1$ is closed with the relative topology of $H_r$. Particularly, $(u_0,v_0) \in X_1$, thus $X_1$ is non-empty and completed.

Next we claim that for $\forall \delta >0$, there exists  some $(f_\delta, g_\delta) \in X_1$, such that
\be \lab{308}
\varphi(f, g)>\varphi(f_\delta,g_\delta)-\frac{\varepsilon}{\delta} \|(f,g)-(f_\delta, g_\delta)\|_H,   \quad \;\forall (f,g)\in X_1\backslash (f_\delta,g_\delta).
\ee
If not, then there exists $\delta_0>0$, such that for any $(f_{\delta_{0}},g_{\delta_0}) \in X_1$, we could find some $(f,g) \in X_1\backslash (f_{\delta_{0}},g_{\delta_0})$, such that
\be \lab{309}
\varphi(f,g)\leqslant  \varphi(f_{\delta_0}, g_{\delta_0})-\frac{\varepsilon}{\delta}\|(f,g)-(f_{\delta_0}, g_{\delta_0})\|_H.
\ee
This enables  us  to define a map $h: X_1 \to X_1$
\be \lab{310}
h(f_{\delta_0}, g_{\delta_0})= (f,g) \neq (f_{\delta_0}, g_{\delta_0}),
\ee
 satisfying
\be  \lab{311}
\frac{\varepsilon}{\delta_0}\|h(f_{\delta_0}, g_{\delta_0})-(f_{\delta_0}, g_{\delta_0})\|_H\leqslant \varphi(f_{\delta_0}, g_{\delta_0})-\varphi \s{h(f_{\delta_0}, g_{\delta_0})}, \quad \forall (f_{\delta_0}, g_{\delta_0})\in X_1.
\ee
 However, by the well known Caristi fixed point theorem, $h$ should have a fixed point, a contradiction.

In the following, for $\forall \delta>0$, let $(f_\delta,g_\delta)$ be the one given by \eqref{308}. Since $(f_\delta,g_\delta) \in X_1$, we have that
\be\lab{312}
\varphi(f_\delta,g_\delta)+\frac{\varepsilon}{\delta}\|(u_0,v_0)-(f_\delta,g_\delta)\|_H\leqslant  \varphi(u_0,v_0).
\ee
Moreover,
\be \lab{313}
\varphi(f_\delta,g_\delta)\leqslant \varphi(u_0,v_0).
\ee
We claim that
\be \lab{314}
\|(u_0,v_0)-(f_\delta,g_\delta)\|_H\leqslant \delta.
\ee
If not, we have
\be \lab{315}
\begin{split}
\varphi(u_0,v_0)\geqslant & \varphi(f_\delta,g_\delta)+\frac{\varepsilon}{\delta}\|(u_0,v_0)-(f_\delta,g_\delta)\|_H\\
>&\varphi(f_\delta,g_\delta)+\varepsilon\\
\geqslant &\inf_{(f,g)\in S(a,b)}\varphi(f,g)+\varepsilon,
\end{split}
\ee
a contradiction to the choice of $(u_0,v_0)$. Next, we shall prove that
\be \lab{316}
\varphi(f,g)>\varphi(f_\delta,g_\delta)-\frac{\varepsilon}{\delta}\|(f,g)-(f_\delta,g_\delta)\|_H,\; \quad \forall (f,g)\in S_r(a,b)\backslash (f_\delta,g_\delta).
\ee
By the definition of $(f_\delta,g_\delta)$, see \eqref{308},  we only need to discuss the cases of $(f,g) \in S_r(a,b) \setminus X_1$. Firstly,  if $(f,g)\in S_r(a,b)\backslash X_1$ with $f,g \geqslant 0$,  then
\be \lab{317}
\varphi(f,g)+\frac{\varepsilon}{\delta}\|(u_0,v_0)-(f,g)\|_H > \varphi(u_0,v_0),
\ee
and thus we have that
\be \lab{318}
\begin{split}
\varphi(f,g)>&\varphi(u_0,v_0)-\frac{\varepsilon}{\delta}\|(u_0,v_0)-(f,g)\|_H\\
(\hbox{since}\;(f_\delta,g_\delta)\in X_1)\geqslant & \varphi(f_\delta,g_\delta)+\frac{\varepsilon}{\delta}\|(u_0,v_0)-(f_\delta,g_\delta)\|_H
-\frac{\varepsilon}{\delta}\|(u_0,v_0)-(f,g)\|_H\\
\geqslant  &\varphi(f_\delta,g_\delta)-\frac{\varepsilon}{\delta}\|(f,g)-(f_\delta,g_\delta)\|_H.
\end{split}
\ee
\eqref{316} still holds. On the other hand, if  the non-negativity of $f,g$ is not satisfied, we note that
\be \lab{319}
\varphi(f,g)\geqslant  \varphi(|f|, |g|)
\ee
and
\be \lab{320}
\|(f,g)-(f_\delta,g_\delta)\|_H \geqslant  \|(|f|,|g|)-(f_\delta,g_\delta)\|_H
\ee
since $f_\delta,g_\delta \geqslant 0$. Hence we have that
\be  \lab{321}
\begin{split}
\varphi(f,g)+\frac{\varepsilon}{\delta}\|(f,g)-(f_\delta,g_\delta)\|_H &\geqslant  \varphi(|f|,|g|)+\frac{\varepsilon}{\delta}\|(|f|,|g|)-(f_\delta,g_\delta)\|_H\\
&> \varphi(f_\delta,g_\delta).
\end{split}
\ee
Thus \eqref{316} holds and we  finish the proof of Lemma \ref{l301} by taking $(u_\delta, v_\delta)=(f_\delta,g_\delta)$.

\ep

We note that $\inf\limits_{(f,g) \in S_r(a,b)} \varphi(f,g)=m_r(a,b)=m(a,b) = \inf\limits_{ (u,v) \in {\cal P} \cap  S_r(a,b)} J(u,v) $, hence we could find a positive radial P-S sequence for $J (u,v)$ as following.
\bc \lab{c301}
There exists a radial P-S sequence $\{ (u_n,v_n)    \} \subset S_r(a,b)$, $u_n \geqslant 0, v_n \geqslant 0$, such that
\be \lab{322}
J(u_n,v_n) \to m(a,b)
\ee
and
\be \lab{323}
J'\big|_{S_r(a,b)}(u_n,v_n) \to 0
\ee
as  $n \to \infty$.     In particular, we can require that
\be  \lab{324}
(u_n,v_n) \in {\cal P}.
\ee

\ec

\bp
 Let $ \varepsilon_n \to 0$ and put $\delta_n = \sqrt{\varepsilon_n}$, by Lemma \ref{l301}, there exists a sequence $\{ (\tilde{u}_n,\tilde{v}_n)    \} \subset S_r(a,b)$, $\tilde{u}_n \geqslant 0, \tilde{v}_n \geqslant 0$, which satisfies $\varphi (\tilde{u}_n,\tilde{v}_n) \to m(a,b)$ and $\varphi ' (\tilde{u}_n,\tilde{v}_n)  \to 0 $ as $n \to \infty$. By the definition of $\varphi (u,v)$, for any $(u,v) \in S_r(a,b)$, there exists a  unique $t_{(u,v)}> 0$ such that
 \be \lab{328}
\varphi (u,v) =\max\limits_{t>0}\tilde{J}(t,u,v) =\tilde{J}(t_{(u,v)},u,v).
 \ee
Hence we could choose $t_n :=t_{(\tilde{u}_n,\tilde{v}_n)}$ satisfies
 \be  \lab{325}
J( t_n \star \s{ \tilde{u}_n,\tilde{v}_n})=\max\limits_{t>0}J\s{t \star \s{\tilde{u}_n,\tilde{v}_n}} = \varphi (\tilde{u}_n,\tilde{v}_n).
 \ee
In particular, $P(t_n \star \s{ \tilde{u}_n,\tilde{v}_n}) =0$ and $\partial_{(u,v)}t(u,v) \big|_{(u,v)=(\tilde{u}_n,\tilde{v}_n )} $ exists by the implicit function theorem.

Set $(u_n,v_n)=(t_n \star \s{ \tilde{u}_n,\tilde{v}_n})$, then $(u_n,v_n) \in S_r(a,b)$ and  we have that $J(u_n,v_n)=\varphi (\tilde{u}_n,\tilde{v}_n) \to m(a,b)$.
Since
\be \lab{326}
\frac{\partial}{\partial t}\tilde{J}(t, \tilde{u}_n,\tilde{v}_n)\big|_{t=t_n}=\frac{1}{t}P\big(t \star (\tilde{u}_n,\tilde{v}_n)\big)\big|_{t=t_n}=\frac{1}{t_n}P(t_n \star \s{ \tilde{u}_n,\tilde{v}_n}) =0,
\ee
there holds
\be   \lab{327}
\begin{split}
\varphi '(\tilde{u}_n,\tilde{v}_n ) &=\partial_t\tilde{J}( t, \tilde{u}_n,\tilde{v}_n ) \big|_{t=t_n} \cdot \partial_{(u,v)}t(u,v) \big|_{(u,v)=(\tilde{u}_n,\tilde{v}_n )}  +\partial_{(u,v)}\tilde{J}( t_n, u,v ) \big|_{(u,v)=(\tilde{u}_n,\tilde{v}_n )}  \\
&=0+ \partial_{(u,v)}J( t_n \star (u,v) ) \big|_{(u,v)=(\tilde{u}_n,\tilde{v}_n )} \\
&=J' |_{S_r(a,b)} ( t_n \star (u,v) ) \big|_{(u,v)=(\tilde{u}_n,\tilde{v}_n )}\\
&=J' |_{S_r(a,b)}(u_n,v_n).
\end{split}
\ee
Thus $J'\big|_{S_r(a,b)}(u_n,v_n) \to 0$.

\ep

%-----------------------------------------------------------------------------------------------------------------------------------------------------------------
\vskip0.26in
\section{Boundedness of the Palais-Smale sequence and Lagrange multipliers}\lab{sec:boundedness}
\renewcommand{\theequation}{4.\arabic{equation}}
\renewcommand{\theremark}{4.\arabic{remark}}
\renewcommand{\thedefinition}{4.\arabic{definition}}

By Corollary \ref{c301}, we could find a radial nonnegative P-S sequence $\{(u_n,v_n)   \}$ of $J(u,v)$ on the level $m(a,b)$.
The conclusion \eqref{323} implies that  there exists two sequences of real numbers $(\lambda_{1,n})$ and $(\lambda_{2,n})$ such that
\be  \lab{400.1}
\begin{split}
o(1)||(\varphi, \psi)||_{H}  &= \int_{\RN}  \nabla u_n \cdot \nabla \varphi + \nabla v_n \cdot \nabla \psi  dx  -\mu_1\int_{\RN} |u_n|^{p-2}u_n \varphi dx-\mu_2\int_{\RN} |v_n|^{q-2}v_n \psi dx\\
&\quad -\beta \int_{\RN}\big( u_n \psi +v_n \varphi \big) dx  + \lambda_{1,n} \int_{\RN} u_n \varphi dx + \lambda_{2,n} \int_{\RN} v_n \psi dx
\end{split}
\ee
for every $(\varphi, \psi) \in H$, with $o(1) \to 0$ as $n \to \infty$. Here we show that the P-S sequence and Lagrange multipliers are bounded.

\bl \lab{l401}
Let $\{  (u_n,v_n) \}$ be a radial nonnegative P-S sequence given by Corollary \ref{c301}, then we have  that $\{  (u_n,v_n) \}$ is bounded in $H$. Moreover, the Lagrange multipliers $(\lambda_{1,n})$ and $(\lambda_{2,n})$ in \eqref{400.1} are also bounded.
\el

\bp
By Corollary \ref{c301}, for the  P-S sequence  $\{  (u_n,v_n) \} \subset S_r(a,b)$, there holds
\be \lab{401}
\begin{split}
m(a,b)&=J(u_n,v_n) +o(1) \\
&=\frac{1}{2}\int_{\RN}| \nabla u_n|^2 +|\nabla v_n|^2 dx -\frac{\mu_1}{p} \int_{\RN}|u_n|^p dx -\frac{\mu_2}{q} \int_{\RN}|v_n|^q dx -\beta \int_{\RN}u_nv_n dx +o(1)
\end{split}
\ee
and
\be \lab{402}
P(u_n,v_n) =\int_{\RN}| \nabla u_n|^2 +|\nabla v_n|^2 dx-\mu_1 \gamma_p \int_{\RN}|u_n|^p dx -\mu_2 \gamma_q \int_{\RN}|v_n|^q dx =0.
\ee
By \eqref{402} it implies that
\be \lab{420.3}
\begin{split}
\int_{\RN}| \nabla u_n|^2 +|\nabla v_n|^2 dx &= \mu_1 \gamma_p \int_{\RN}|u_n|^p dx + \mu_2 \gamma_q \int_{\RN}|v_n|^q dx \\
&\geqslant \frac{N}{2} \s{\min \{ p,q  \}-2   } \s{ \frac{\mu_1}{p} \int_{\RN}|u_n|^p dx +\frac{\mu_2}{q} \int_{\RN}|v_n|^q dx   }.
\end{split}
\ee
Combing with \eqref{401}, we have that
\be  \lab{404}
\begin{split}
-C & \leqslant -m(a,b) - \beta \int_{\RN}u_nv_n dx +o(1)\\
&= \frac{\mu_1}{p} \int_{\RN}|u_n|^p dx +\frac{\mu_2}{q} \int_{\RN}|v_n|^q dx-\frac{1}{2}\int_{\RN}| \nabla u_n|^2 +|\nabla v_n|^2 dx \\
& \leqslant \Big[\frac{2}{N \s{\min\{p,q \}-2   }}-\frac{1}{2} \Big]\int_{\RN}| \nabla u_n|^2 +|\nabla v_n|^2 dx,
\end{split}
\ee
thus $| \nabla u_n|_2^2 +|\nabla v_n|_2^2 \leqslant C $ due to the fact of $p,q > 2+ \frac{4}{N}$. Hence,  $\{  (u_n,v_n) \}$ is bounded in $H$.

\vskip0.12in
Next, we estimate the Lagrange multipliers. By taking $(u_n,0)$ as test function for \eqref{400.1}, we obtain that
\be \lab{406}
o(1)\|u_n\|_{H^1}= \int_{\RN}|\nabla u_n|^2dx -\mu_1\int_{\RN}|u_n|^p dx -\beta \int_{\R^N}u_nv_n + \lambda_{1,n} a.
\ee
So
\be \lab{407}
 |\lambda_{1,n}|=\frac{1}{a}\left|o(1)\|u_n\|_{H^1}-\int_{\RN}|\nabla u_n|^2dx+\mu_1\int_{\RN}|u_n|^p dx+\beta\int_{\R^N}u_nv_n\right|<+\infty.
\ee
Similarly, we can prove that $\lambda_{2,n}$ is bounded. Thus the Lagrange multipliers $(\lambda_{1,n})$ and $(\lambda_{2,n})$  are also bounded.

\ep

\vskip0.26in
\section{Compactness of the Palais-Smale sequence}\lab{sec:compactness}
\renewcommand{\theequation}{5.\arabic{equation}}
\renewcommand{\theremark}{5.\arabic{remark}}
\renewcommand{\thedefinition}{5.\arabic{definition}}
In this section, we study the compactness of the Palais-Smale sequence given by Corollary \ref{c301}.
\subsection{The case for $2+\frac{4}{N}<p,q< 2^*$}
\bl \lab{l402}
Assume $2\leqslant N\leqslant 4$ and  $2+\frac{4}{N}<p,q< 2^*$.
Let $\{(u_n,v_n)\}$ be a minimizing P-S sequence given by Corollary 3.1, then up to a subsequence, $(u_n,v_n) \rightarrow (\bar{u},\bar{v})\in S_r(a,b)$  strongly in $H$.
\el

\bp
Note that $(u_n,v_n)\in S_r(a,b)$ is bounded, up to a subsequence, we assume that $(u_n,v_n)\rightharpoonup (\bar{u},\bar{v})\in H_r$. For $2 \leqslant N \leqslant 4$, we have that
 \be \lab{408}
(u_n,v_n)\rightarrow (\bar{u},\bar{v})\; \ \text{in} \;L^p(\R^N)\times L^q(\R^N)
\ee
via the compact embedding of  the radial space.

%For $N=1$, compactness still holds since $\{(u_n,v_n)\}$ are  radially decreasing functions (see e.g. \cite[ Proposition 1.7.1]{CT2003}).

Firstly we show that $(\bar{u},\bar{v}) \neq (0,0)$. If not, by $(u_n, v_n)\in \mathcal{P}$, we have that
\be  \lab{409}
\int_{\RN}| \nabla u_n|^2 +|\nabla v_n|^2 dx=\mu_1 \gamma_p \int_{\RN}|u_n|^p dx + \mu_2 \gamma_q \int_{\RN}|v_n|^q dx \to 0,
\ee
a contradiction to Lemma \ref{lemma:bounded-away}.

Up to a subsequence, we assume that $\lambda_{i,n}\rightarrow \lambda_i, i=1,2.$ Then $ H_r \ni(\bar{u},\bar{v})$ is a nonnegative solution to
\be \lab{412}
\begin{cases}
-\Delta u+\lambda_1 u=\mu_1 |u|^{p-2}u+\beta v,\;\hbox{in}\;\R^N,\\
-\Delta v+\lambda_2 v=\mu_2 |v|^{q-2}v+\beta u,\;\hbox{in}\;\R^N.
\end{cases}
\ee
So if $\bar{u}=0$, we also have that $\bar{v}=0$ due to the fact $\beta>0$, a contradiction to $(\bar{u},\bar{v})\neq (0,0)$. Hence by the strong maximum principle we have that $\bar{u}>0, \bar{v}>0$.

Furthermore, we have $\lambda_i>0$ for $ i=1,2$. Indeed, if $\lambda_1 \leqslant 0$, then  \eqref{412} implies that
\be \lab{413}
-\Delta u\geqslant 0\;\hbox{in}\;\R^N, u\in L^2(\R^N),
\ee
a contradiction to \cite[Lemma A.2]{ikoma2014compactness}. Hence, $\lambda_1>0$. Similarly, we have $\lambda_2>0$.

Now we prove the compactness of $\{(u_n,v_n)\}$.  Combine $(u_n,v_n) \in  { \cal P}$ with $\eqref{408}$, we have that
\be \lab{414.1}
\begin{split}
\int_{\RN}| \nabla u_n|^2 +|\nabla v_n|^2 dx&=\mu_1 \gamma_p \int_{\RN}|u_n|^p dx +\mu_2 \gamma_q \int_{\RN}|v_n|^q dx\\
&=\mu_1 \gamma_p \int_{\RN}|\bar{u}|^p dx +\mu_2 \gamma_q \int_{\RN}|\bar{v}|^q dx+o(1)
\end{split}
\ee
Observe that $(\bar{u},\bar{v})$ is a solution to \eqref{412}, we also have $(\bar{u},\bar{v})\in \mathcal{P}$ and then
\be \lab{415}
\int_{\RN}| \nabla u_n|^2 +|\nabla v_n|^2 dx \to \int_{\RN}| \nabla \bar{u}|^2 +|\nabla \bar{v}|^2 dx  \ \ \text{as} \ \ n \to \infty.
\ee
Hence,
$$(u_n,v_n)\rightarrow (\bar{u},\bar{v})\; \text{in} \;D_{0}^{1,2}(\R^N)\times D_{0}^{1,2}(\R^N).$$
Set $|\bar{u}|_2^2=\bar{a}\in (0,a], |\bar{v}|_2^2=\bar{b}\in (0,b] $. By
\be \lab{416}
\begin{cases}
|\nabla u_n|_2^2+\lambda_{1,n}|u_n|_2^2=\mu_1|u_n|_p^p+\beta \langle u_n,v_n\rangle,\\
|\nabla v_n|_2^2+\lambda_{2,n}|v_n|_2^2=\mu_2|v_n|_q^q+\beta \langle u_n,v_n\rangle,\\
|\nabla \bar{u}|_2^2+\lambda_1|\bar{u}|_2^2=\mu_1|\bar{u}|_p^p+\beta \langle \bar{u} ,\bar{v}\rangle,\\
|\nabla \bar{v}|_2^2+\lambda_2|\bar{v}|_2^2=\mu_2|\bar{v}|_q^q+\beta \langle \bar{u},\bar{v}\rangle,\\
|\nabla u_n|_2^2\rightarrow |\nabla \bar{u}|_2^2, |\nabla v_n|_2^2\rightarrow |\nabla \bar{v}|_2^2,\\
|u_n|_p^p\rightarrow |\bar{u}|_p^p, |v_n|_q^q\rightarrow |\bar{v}|_q^q,
\lambda_{i,n}\rightarrow \lambda_i>0, i=1,2,
\end{cases}
\ee
we have that
\be \lab{417}
\lambda_1(a-\bar{a})=\beta \lim_{n\rightarrow \infty} \langle u_n-u,v_n-v\rangle =\lambda_2(b-\bar{b})=:\kappa.
\ee
By Fatou's lemma, we have that $\kappa\geqslant  0$.

{\bf We claim that $\kappa=0$.}
If not,  $\kappa>0$ and thus
\be \lab{418}
\bar{a}<a, \bar{b}<b, \lambda_1>0, \lambda_2>0.
\ee
And then by Schwartz inequality,
\be \lab{419}
\begin{split}
2\sqrt{\lambda_1\lambda_2} \sqrt{(a-\bar{a})(b-\bar{b})}\leqslant&\lambda_1(a-\bar{a})+\lambda_2(b-\bar{b})\\
=&2\beta \lim_{n\rightarrow \infty} \langle u_n-\bar{u},v_n-\bar{v}\rangle\\
\leqslant & 2\beta \sqrt{(a-\bar{a})(b-\bar{b})},
\end{split}
\ee
which implies that
\be \lab{420}
\beta\geqslant \sqrt{\lambda_1\lambda_2}.
\ee
Recalling \eqref{412}, we see that $(\tilde{u}, \tilde{v}):=(\sqrt{\lambda_2}\bar{u}, \sqrt{\lambda_1}\bar{v})$ are positive radial functions of $H_r$
satisfying
\be \lab{421}
\begin{cases}
-\Delta  \tilde{u}+\lambda_1\tilde{u}\geqslant \lambda_2 \tilde{v},\\
-\Delta \tilde{v}+\lambda_2 \tilde{v}\geqslant \lambda_1 \tilde{u},
\end{cases}
\ee
and thus
\be \lab{422}
-\Delta (\tilde{u}+\tilde{v})\geqslant  0, (\tilde{u}+\tilde{v})\in L^2(\R^N),
\ee
a contradiction to \cite[Lemma A.2]{ikoma2014compactness} again. Hence, the claim about  $\kappa=0$ is proved.

Recalling that $\lambda_1>0$ and

\be \lab{423}
0=\kappa=\lambda_1(a-\bar{a}),
\ee
we obtain that $\bar{a}=a$. Similarly, we have that $\bar{b}=b$.
That is $(u_n,v_n)\rightarrow (\bar{u},\bar{v})$ in $H$.
\ep

%-----------------------------------------------------------------------------------------------------------------------------------------------------------------
\subsection{The case for $2+\frac{4}{N}<p< 2^*, q=2^*$}

In this subsection, we consider the case of $2+\frac{4}{N}<p<2^*,  q=2^*$ for $N=3,4$. By Corollary \ref{c301} and Lemma \ref{l401}, we could also find a bounded nonnegative radial P-S sequence  with the Lagrange multipliers are bounded. Since  $(u_n,v_n)\in S_r(a,b)$ is bounded, up to a subsequence, we assume that $(u_n,v_n)\rightharpoonup (\bar{u},\bar{v})\in H_r$. Moreover, we have that
$u_n\rightarrow \bar{u}\; \ \text{in} \;L^p(\R^N)$
by the compact embedding for the radial space. In addition, up to a subsequence, we may assume that $\lambda_{i,n}\rightarrow \lambda_i, i=1,2.$ Then $H_r\ni(\bar{u},\bar{v})$ is a solution to
\be \lab{4201}
\begin{cases}
-\Delta u+\lambda_1 u=\mu_1 |u|^{p-2}u+\beta v,\;\hbox{in}\;\R^N,\\
-\Delta v+\lambda_2 v=\mu_2 |v|^{q-2}v+\beta u,\;\hbox{in}\;\R^N.
\end{cases}
\ee
Since $\beta > 0$, applying a similar argument as the proof of Lemma \ref{l402}, we have that one of the following holds:
\begin{itemize}
\item[(i)]${\bf (\bar{u},\bar{v})=(0,0)}$;
\item[(ii)]${\bf \bar{u}>0,\bar{v}>0}$.
\end{itemize}

\bl  \lab{l403}
Assume that  $a>0,b>0$ satisfyies
\be \lab{4203}
m(a,b)< \frac{1}{N}\frac{{ \cal S}^{\frac{N}{2}}}{\mu_2^{\frac{N}{2}-1}} -\beta \sqrt{ab},
\ee
then ${\bf \bar{u}>0,\bar{v}>0}$.
\el
\bp
Suppose that $(\bar{u},\bar{v}) = (0,0)$, then $u_n \to 0$ in $L^p(\RN)$. Since $(u_n,v_n) \in {\cal P}$, we have that
\be \lab{4204}
\int_{\RN}| \nabla u_n|^2 +|\nabla v_n|^2 dx  -\mu_2  \int_{\RN}|v_n|^{2^*} dx =o(1).
\ee
Notice that  $\{ v_n\}$ is bounded in $H^1(\RN)$, up to a subsequence, we could assume  that $|v_n|_{2^*} \to \ell \in \R $. By \eqref{4204} and Sobolev inequality, there holds $  {\cal S} \ell^2 \leqslant \mu_2 \ell^{2^*}$, hence we deduce that
$$ \text{either} \  \ell=0, \ \text{or} \ \ell \geqslant \s{\frac{ {\cal S }}{\mu_2} }^{\frac{N-2}{4}}. $$
Let us suppose firstly that $\ell=0$. By $(u_n, v_n)\in \mathcal{P}$, we have that
\be  \lab{4205}
\int_{\RN}| \nabla u_n|^2 +|\nabla v_n|^2 dx=\mu_1 \gamma_p \int_{\RN}|u_n|^p dx + \mu_2  \int_{\RN}|v_n|^{2^*} dx \to 0,
\ee
 a contradiction to Lemma \ref{lemma:bounded-away}.

Secondly, if  $\ell \geqslant \s{\frac{ {\cal S }}{\mu_2} }^{\frac{N-2}{4}}$, by  \eqref{4204}, the following holds
\be \lab{4206}
\begin{split}
m(a,b)&=J(u_n,v_n)+o(1)\\
&=\frac{1}{2}\int_{\RN}| \nabla u_n|^2 +|\nabla v_n|^2 dx  -\frac{\mu_2}{2^*} \int_{\RN}|v_n|^{2^*} dx -\beta \int_{\RN}u_nv_n dx +o(1)\\
&=\frac{\mu_2}{N}\int_{\RN}|v_n|^{2^*} dx -\beta \int_{\RN}u_nv_n dx +o(1)\\
& \geqslant \frac{1}{N}\frac{{ \cal S}^{\frac{N}{2}}}{\mu_2^{\frac{N}{2}-1}} -\beta \sqrt{ab},
\end{split}
\ee
also a contradiction.
\ep

 Now we have that $\bar{u}, \bar{v}>0$ and we could prove the convergence of $\{(u_n,v_n) \}$.
\bl \lab{l404}
Let $\{(u_n,v_n)\}$ be a minimizing P-S sequence given by Corollary \ref{c301} with $$m(a,b)+\beta \sqrt{ab}< \frac{1}{N}\frac{{ \cal S}^{\frac{N}{2}}}{\mu_2^{\frac{N}{2}-1}},$$
 then up to a subsequence, $(u_n,v_n) $   converges  to some   $ (\bar{u},\bar{v})\in S(a,b)$ strongly in $H$.
\el
\bp
By Lemma \ref{l403}, $(\bar{u},\bar{v}) \neq (0,0)$. Set $|\bar{u}|_2^2=\bar{a}\in (0,a], |\bar{v}|_2^2=\bar{b}\in (0,b] $. Applying a similar arguments as the proof of Lemma \ref{l402}, one can prove  that $(\bar{u}, \bar{v})$ is a  positive solution to
\be \lab{4211}
\begin{cases}
-\Delta u+\lambda_1 u=\mu_1 |u|^{p-2}u+\beta v,\;\hbox{in}\;\R^N,\\
-\Delta v+\lambda_2 v=\mu_2 |v|^{2^*-2}v+\beta u,\;\hbox{in}\;\R^N
\end{cases}
\ee
 with $\lambda_i>0,i=1,2$. Put $(\tilde{u}_n,\tilde{v}_n)=(u_n -\bar{u},v_n- \bar{v})$, by $P(u_n,v_n) =P(\bar{u},\bar{v})=0$, we have
\be \lab{4212}
\int_{\RN}| \nabla \tilde{u}_n|^2 +|\nabla \tilde{v}_n|^2 dx  -\mu_2  \int_{\RN}|\tilde{v}_n|^{2^*} dx =o(1).
\ee
Furthermore,
\be \lab{4213}
\int_{\RN}| \nabla \tilde{u}_n|^2dx + \lambda_1 (a-\bar{a}) = \beta \int_{\RN}\tilde{u}_n\tilde{v}_n dx +o(1),
\ee
\be\lab{4214}
\int_{\RN}| \nabla \tilde{v}_n|^2dx + \lambda_2 (b-\bar{b}) =  \mu_2 \int_{\RN} |\tilde{v}_n|^{2^*} dx +     \beta \int_{\RN}\tilde{u}_n\tilde{v}_n dx +o(1).
\ee
Notice that  $\{ \tilde{ v}_n\}$ is bounded in $H^1(\RN)$, up to a subsequence, we can assume that $|\tilde{v}_n|_{2^*} \to \ell \in \R $. By \eqref{4212} and Sobolev inequality, there holds $  {\cal S} \ell^2 \leqslant \mu_2 \ell^{2^*}$, hence we deduce that
$$ \text{either} \  \ell=0, \ \text{or} \ \ell \geqslant \s{\frac{ {\cal S }}{\mu_2} }^{\frac{N-2}{4}}. $$
Let us suppose firstly that $\ell \geqslant \s{\frac{ {\cal S }}{\mu_2} }^{\frac{N-2}{4}}$, then
\be \lab{513}
\begin{split}
m(a,b)&=J(u_n,v_n)+o(1)\\
&=\frac{1}{2}\int_{\RN}| \nabla \bar{u}|^2 +|\nabla \bar{v}|^2 dx  -\frac{\mu_1}{p} \int_{\RN}|\bar{u}|^{p} dx -\frac{\mu_2}{2^*} \int_{\RN}|\bar{v}|^{2^*} dx \\
& \quad \ + \frac{1}{2}\int_{\RN}| \nabla \tilde{u}_n|^2 +|\nabla \tilde{v}_n|^2 dx  -\frac{\mu_2}{2^*} \int_{\RN}|\tilde{v}_n|^{2^*} dx -\beta \int_{\RN}u_nv_n dx +o(1)\\
& \geqslant \frac{\mu_1}{2}(\gamma_p-\frac{2}{p})\int_{\RN}|\bar{u}|^{p} dx + \frac{\mu_2}{N}\int_{\RN}|\bar{v}|^{2^*} dx +\frac{1}{N}\frac{{ \cal S}^{\frac{N}{2}}}{\mu_2^{\frac{N}{2}-1}}-\beta\sqrt{a b}\\
&>\frac{1}{N}\frac{{ \cal S}^{\frac{N}{2}}}{\mu_2^{\frac{N}{2}-1}}-\beta\sqrt{ab},
\end{split}
\ee
a contradiction. Hence, $\ell= 0 $, and thus by \eqref{4212}, we have
$$\int_{\RN}| \nabla \tilde{u}_n|^2 +|\nabla \tilde{v}_n|^2 dx \rightarrow 0\; \text{as} \;n\rightarrow \infty,$$
which implies that
$$u_n\rightarrow \bar{u}, \bar{v}_n\rightarrow \bar{v}\;\hbox{strongly in
$D_{0}^{1,2}(\R^N)$}.$$
Going to a subsequence, we let $\displaystyle \kappa:=\lim_{n\rightarrow \infty}\beta \int_{\RN}\tilde{u}_n\tilde{v}_n dx \geqslant 0$. By \eqref{4213} and \eqref{4214}, we have that
\be \lab{4215}
\lambda_1 (a-\bar{a})=\lambda_2 (b-\bar{b})=\kappa.
\ee
So we can apply a similar argument as the proof of Lemma \ref{l402}, and obtain that
\be \lab{4216}
\kappa=0\;\hbox{and}\;(\bar{a}, \bar{b})=(a,b).
\ee
Consequently, we have that $(u_n,v_n) \to (\bar{u},\bar{v})$ strongly in $H$.
\ep

%%%%%%%%%%%%%%%%%%%%%%%%%%%%%%%%%%%%%%%%%%%%%%%%%%%%%%%%%%%%%%%%%%%%%%%%%%%%%%%%%%
\vskip0.26in
\section{Proof of Theorems}\lab{sec:proofs}
\renewcommand{\theequation}{6.\arabic{equation}}
\renewcommand{\theremark}{6.\arabic{remark}}
\renewcommand{\thedefinition}{6.\arabic{definition}}

%%%%%%%%%%%%%%%%%%%%%%%%%%%%%%%%%%%%%%%%%%%%%%%%%%%%%%%%%%%%%%%%%%%%%%%%%%%%%%%%%%

\begin{proof}[Proof of Theorem \ref{thm1}]
By Corollary \ref{c301}, there exists a nonnegative radial minimizing  P-S sequence $\{ (u_n,v_n)\} \subset S_r(a,b)\cap {\cal P}$, such that
$$
J(u_n,v_n) \to m(a,b)\;\hbox{and}\; J'\big|_{S_r(a,b)}(u_n,v_n) \to 0.
$$
By Lemma \ref{l402}, $(u_n,v_n)\rightarrow (\bar{u},\bar{v})$ strongly in $H$, which is a solution to \eqref{101} with $|\bar{u}|_2^2=a, |\bar{v}|_2^2=b,  \lambda_1,\lambda_2 >0$. In particular, $J(\bar{u},\bar{v}) =m(a,b)$, hence $( \bar{u},\bar{v})$ is a ground state solution in the sense that
\be
\begin{split}
J(\bar{u},\bar{v}) &=\inf \{ (u,v): (u,v) \in {\cal P} \cap S(a,b)    \} \\
& =\inf \{ (u,v): (u,v) \ \  \text{is a solution of \eqref{101}-\eqref{102} for some }\lambda_1,\lambda_2   \}
\end{split}
\ee
holds.
\end{proof}

\begin{proof}[Proof of Theorem \ref{thm2}]
By Corollary \ref{c301},there exists a nonnegative radial minimizing  P-S sequence $\{ (u_n,v_n)\} \subset S_r(a,b)\cap {\cal P}$, such that
$$
J(u_n,v_n) \to m(a,b)\;\hbox{and}\; J'\big|_{S_r(a,b)}(u_n,v_n) \to 0.
$$
If  $$m(a,b)< \frac{1}{N}\frac{{ \cal S}^{\frac{N}{2}}}{\mu_2^{\frac{N}{2}-1}} -\beta \sqrt{ab},$$ then by Lemma \ref{l403} and \ref{l404}, $(u_n,v_n)\rightarrow (\bar{u},\bar{v})$ strongly in $H$, which is a solution to \eqref{101} with $|\bar{u}|_2^2=a, |\bar{v}|_2^2=b,  \lambda_1,\lambda_2 >0$. In particular, $J(\bar{u},\bar{v}) =m(a,b)$, hence $( \bar{u},\bar{v})$ is a ground state solution in the sense that
\be
\begin{split}
J(\bar{u},\bar{v}) &=\inf \{ (u,v): (u,v) \in {\cal P} \cap S(a,b)    \} \\
& =\inf \{ (u,v): (u,v) \ \  \text{is a solution of \eqref{101}-\eqref{102} for some }\lambda_1,\lambda_2   \}
\end{split}
\ee
holds.
\end{proof}

\begin{proof}[Proof of Corollary \ref{cro-th2}]
We claim that
\be \lab{4206.1}
 m(a,b) \leqslant m(a,0)=\Big( \frac{1}{2}-\frac{1}{p\gamma_p} \Big) \Big( \gamma_p  {\cal C}_{N,p}\mu_1 a^\frac{p-p\gamma_p}{2} \Big)^\frac{2}{2-p\gamma_p}.
\ee
  Indeed, let $u =u_{\mu_1,p,a}$ be the unique positive solution to \eqref{202} with $ v \in S(b)$. Then $(u, s \star v)  \in S(a,b)$ for any $s \in \R$. Let $t_s := t_{(u, s \star v)}$ be the unique positive number such that $P( t_{(u, s \star v)} \star (u, s \star v)) \equiv 0$. Then we have that
\be  \lab{4207}
\begin{split}
0&=P( t_s \star (u, s \star v)) \\
&=t_s^2\int_{\RN}| \nabla u|^2 dx + t_s^2 \cdot s^2 \int_{\RN}| \nabla v|^2 dx -\mu_1 \gamma_p t_s^{p \gamma_p}\int_{\RN}| u|^p dx- \mu_2 t_s^{2^*} \cdot s^{2^*} \int_{\RN}| \nabla v|^2 dx,
\end{split}
\ee
which means that
\be  \lab{4208}
\int_{\RN}| \nabla u|^2 dx +  s^2 \int_{\RN}| \nabla v|^2 dx \geqslant \mu_1 \gamma_p t_s^{p \gamma_p-2}\int_{\RN}| u|^p dx.
\ee
Hence as $s \to 0^+$, $t_s$ is bounded due to the fact $p\gamma_p>2$. Observing that
\be \lab{4209}
\begin{split}
m(a,b) &\leqslant J(t_s \star (u, s \star v)) \\
&\leqslant J_{\mu_1,p}(t_s \star u) + \frac{(t_s \cdot s)^2}{2}\int_{\RN}| \nabla v|^2 dx -\frac{\mu_2}{2^*}(t_s \cdot s)^{2^*}\int_{\RN}| v|^{{2^*}} dx,
\end{split}
\ee
by letting  $s \to 0^+$, we obtain that
$$m(a,b)\leqslant \lim_{s\rightarrow 0^+}J_{\mu_1,p}(t_s \star u)\leqslant J_{\mu_1,p}( u)=m(a,0).$$
Hence the claim is proved. Combing with \eqref{106.2}, we obtain that
$$m(a,b)\leqslant m(a,0)=\Big( \frac{1}{2}-\frac{1}{p\gamma_p} \Big) \Big( \gamma_p  {\cal C}_{N,p}\mu_1 a^\frac{p-p\gamma_p}{2} \Big)^\frac{2}{2-p\gamma_p}<\frac{1}{N}\frac{{ \cal S}^{\frac{N}{2}}}{\mu_2^{\frac{N}{2}-1}} -\beta \sqrt{ab},$$
so the conclusions of Theorem \ref{thm2} hold.
\end{proof}

%------------------------------------------------------------------------------------------------

This critical case is much easier than in the previous two cases. Here  we prove a stronger conclusion without the mass constraint. And Theorem \ref{thm3} follows directly by the following lemma.
\bl \lab{405}
 Assume that $ N=3,4$ and  $\mu_1, \mu_2,\beta>0$, then
\be\lab{4301}
\begin{cases}
-\Delta u + \lambda_1 u= \mu_1|u|^{2^*-2}u+\beta v \quad &\hbox{in}\;\RN, \\
-\Delta v + \lambda_2 v= \mu_2|v|^{2^*-2}v+\beta u \quad &\hbox{in}\;\RN, \\
u \in H^1(\RN), v \in H^1(\RN)
\end{cases}
\ee
has no nontrivial positive solutions.
\el

\bp
If not, suppose that $(\bar{u},\bar{v})$ solves \eqref{4301} with $\int_{\RN}\bar{u}^2 dx=a>0, \int_{\RN}\bar{v}^2dx =b>0$.
Firstly, we show that $\lambda_i>0$ for $i=1,2$. Suppose $\lambda_1\leqslant 0$, then $-\Delta u
\geqslant 0, u\in L^2(\R^N)$ and by \cite[Lemma A.2]{ikoma2014compactness}, there holds $u \equiv0$, which is a contradiction. Next, by the Pohozaev identity and Nehari identity, there holds that
\be  \lab{1102}
 \lambda_1a+ \lambda_2b=  \lambda_1 \int_{\RN}\bar{u}^2dx + \lambda_2 \int_{\RN}\bar{v}^2 dx= 2\beta \int_{\RN}\bar{u}\bar{v} dx \leqslant 2\beta \sqrt{ab},
\ee
hence $\beta \geqslant \frac{1}{2} \s{\lambda_1 \sqrt{\frac{a}{b}} +\lambda_2 \sqrt{\frac{b}{a}}} \geqslant \sqrt{\lambda_1 \lambda_2}$. However, by setting  $(U,V) =(\sqrt{\lambda_2} \bar{u}, \sqrt{\lambda_1}\bar{v})$, then $(U,V)$ satisfies that
 \be  \lab{1103}
\begin{cases}
-\Delta U + \lambda_1 U= \mu_1  \lambda_2^{\frac{2}{N-2}} |U|^{2^*-2}U+\beta \sqrt{\frac{\lambda_2}{\lambda_1}}V \quad &\hbox{in}\;\RN, \\
-\Delta V + \lambda_2 V= \mu_2\lambda_1^{\frac{2}{N-2}} |V|^{2^*-2}V+\beta \sqrt{\frac{\lambda_1}{\lambda_2}}U \quad &\hbox{in}\;\RN, \\
\int_{\RN}U^2 =\lambda_1  a, \;\;\; \int_{\RN}V^2 =\lambda_2  b,
\end{cases}
\ee
and $U>0, V>0$.
It follows that
 $$-\Delta (U+V) \geqslant 0,  U+V\in L^2(\R^N).$$
A contradiction to  \cite[Lemma A.2]{ikoma2014compactness} again.
\ep

\br  We remark that in \cite[Remark 1.4]{chen2012ground},  when assuming that  $\lambda_1>0, \lambda_2>0, \mu_1=\mu_2=1$ and $0<\beta<\sqrt{\lambda_1\lambda_2}$, Chen and Zou proved that  \eqref{4301} has no nontrivial solution  via the  Pohozaev identity.
\er

%%%%%%%%%%%%%%%%%%%%%%%%%%%%%%%%%%%%%%%%%%%%%%%%%%%%%%%%%%%%%%%%%%%%%%%%%%%%%%%%%%%%%%%%%%%%%%%%%%%%%%%%%%%%%%%%%%%%%%%%%%%

 \end{document}